\documentclass{article}

\usepackage{arxiv}

\usepackage[utf8]{inputenc} 
\usepackage[T1]{fontenc}    
\usepackage[hypertexnames=true]{hyperref}       
\usepackage{url}            
\usepackage{booktabs}       
\usepackage{amsmath, amsfonts}       
\usepackage{amsthm}
\usepackage{nicefrac}       
\usepackage{microtype}      
\usepackage{cleveref}       
\usepackage{graphicx}
\usepackage{doi}

\usepackage{amssymb}
\usepackage{algorithmic}
\usepackage[ruled,linesnumbered,vlined]{algorithm2e}
\usepackage{array}
\usepackage[caption=false,font=normalsize,labelfont=sf,textfont=sf]{subfig}
\usepackage{textcomp}
\usepackage{stfloats}
\usepackage{verbatim}
\usepackage{cite}
\usepackage{multirow} 
\usepackage{scalerel}
\usepackage{tikz}
\usepackage[title,titletoc]{appendix}
\usepackage{multirow}
\usepackage{amssymb}
\usepackage{pifont}
\usepackage{caption}

\DeclareMathOperator*{\argmin}{\arg\min}

\theoremstyle{plain}

\theoremstyle{definition}

\theoremstyle{remark}

\title{Memetic Differential Evolution Methods for Semi-Supervised Clustering\thanks{Under consideration at Pattern Recognition Letters.}}

\author{
	\hspace{1mm}Pierluigi Mansueto \\
	Global Optimization Laboratory (GOL) \\
	Department of Information Engineering \\
	University of Florence \\
	Via di Santa Marta, 3, 50139, Florence, Italy \\
	\texttt{pierluigi.mansueto@unifi.it} \\
	\And
	\hspace{1mm}Fabio Schoen \\
	Global Optimization Laboratory (GOL) \\
	Department of Information Engineering \\
	University of Florence \\
	Via di Santa Marta, 3, 50139, Florence, Italy \\
	\texttt{fabio.schoen@unifi.it} \\
}

\renewcommand{\headeright}{Pierluigi Mansueto and Fabio Schoen}
\renewcommand{\undertitle}{}
\renewcommand{\shorttitle}{}

\hypersetup{
pdftitle={Memetic Differential Evolution Methods for Semi-Supervised Clustering},
pdfsubject={},
pdfauthor={Pierluigi Mansueto, Fabio Schoen},
pdfkeywords={Semi-Supervised Clustering, Minimum Sum-of-Squares Clustering, Memetic Differential Evolution, Global Optimization},
}

\begin{document}
\maketitle

\begin{abstract}
	In this paper, we propose an extension for semi-supervised Minimum Sum-of-Squares Clustering (\texttt{MSSC}) problems of \texttt{MDEClust}, a memetic framework based on the Differential Evolution paradigm for unsupervised clustering. In semi-supervised \texttt{MSSC}, background knowledge is available in the form of (instance-level) ``must-link'' and ``cannot-link'' constraints, each of which indicating if two dataset points should be associated to the same or to a different cluster, respectively. The presence of such constraints makes the problem at least as hard as its unsupervised version and, as a consequence, some framework operations need to be carefully designed to handle this additional complexity: for instance, it is no more true that each point is associated to its nearest cluster center. As far as we know, our new framework, called \texttt{S-MDEClust}, represents the first memetic methodology designed to generate a (hopefully) optimal \textit{feasible} solution for semi-supervised \texttt{MSSC} problems. Results of thorough computational experiments on a set of well-known as well as synthetic datasets show the effectiveness and efficiency of our proposal. 
\end{abstract}

\MSCs{90C11 \and 90C30 \and 90C59}

\section{Introduction}
\label{sec::introduction}

Clustering aims to group together data records based on suitable similarity and dissimilarity criteria, and, over the years, it has been recognized as one of the most relevant statistical techniques (see, e.g., \cite{JAIN2010651}). 

In this paper, we are interested in the \textit{Semi-supervised Euclidean Minimum Sum-of-Squares Clustering} (\texttt{MSSC}) model, where, differently from the (well-known) unsupervised one, some a-priori information on the clusters is available. Successful application of this problem can be found in different domains, e.g., gene \cite{Pensa08} and document \cite{Huang06} clustering. Let us consider $\mathcal{D} = \{p_1, p_2,\ldots, p_N\}$, a dataset of $N$ $d$-dimensional points, and $K$ disjoint sets, called \textit{clusters}, each of which identified by a \textit{center} $y_k \in \mathbb{R}^d$; let us denote the set of the cluster centers as $\mathcal{Y} = \{y_1,\ldots, y_K\}$. A possible formulation for the \texttt{MSSC} problem is thus the following: 
\begin{equation}
	\label{eq::mssc-prob}
	\min_{\substack{x \in \mathcal{X}\cap\mathcal{L}\\y_1,\ldots, y_K \in \mathbb{R}^d}}\;\sum_{i = 1}^{N}\sum_{k = 1}^{K} x_{ik}\|p_i - y_k\|^2,
\end{equation}
where each binary variable $x_{ik}$ indicates whether the point $p_i$ is assigned to the cluster $C_k$,
\begin{equation}
	\label{eq::feasible_set_X}
	\mathcal{X}  = \left\{x \in \{0, 1\}^{N \times K} \middle\vert
	\begin{aligned}
		&\sum_{k = 1}^{K} x_{ik} = 1,\ \forall i \in \{1,\ldots, N\} \\ &\sum_{i = 1}^{N} x_{ik} \ge 1,\ \forall k \in \{1,\ldots, K\}
	\end{aligned}
	\right\},
\end{equation} 
\begin{gather}
	\label{eq::feasible_set_L}
	\mathcal{L} = \left\{x \middle\vert
	\begin{aligned}
		& x_{ik} = x_{jk},\ \forall (i, j) \in \mathcal{M}_\mathcal{L},\ \forall k \in \{1,\ldots, K\} \\
		& x_{ik} + x_{jk} \le 1,\ \forall (i, j) \in \mathcal{C}_\mathcal{L},\ \forall k \in \{1,\ldots, K\}
	\end{aligned}
	\right\}
\end{gather}
and $\|\cdot\|$ indicates the Euclidean norm in $\mathbb{R}^d$. Here, the a-priori information is given in the form of (instance-level) ``must-link'' and ``cannot-link'' constraints in the set $\mathcal{L}$: for each $(i_m, j_m) \in \mathcal{M}_\mathcal{L}$, we require that $p_{i_m}$ and $p_{j_m}$ belong to the same cluster, while, for each $(i_c, j_c) \in \mathcal{C}_\mathcal{L}$, the points $p_{i_c}$ and $p_{j_c}$ should be assigned to different clusters. Clustering is then performed by choosing the centers $y_1,\ldots, y_K$ so that the sum of squared distances from each sample point to its associated center is minimized \cite{hansen_cluster_1997}; each point $p_i$ should be assigned to one and only one center $y_k$ (first set of constraints in \eqref{eq::feasible_set_X}), each cluster should contain at least one sample (second set of constraints) and each instance-level constraints \eqref{eq::feasible_set_L} should be satisfied. From now on, let us assume that problem \eqref{eq::mssc-prob} has always a feasible solution. 

\begin{figure}[h]
	\centering
	\subfloat[Unconstrained \texttt{MSSC}]{\includegraphics[width=0.4\textwidth]{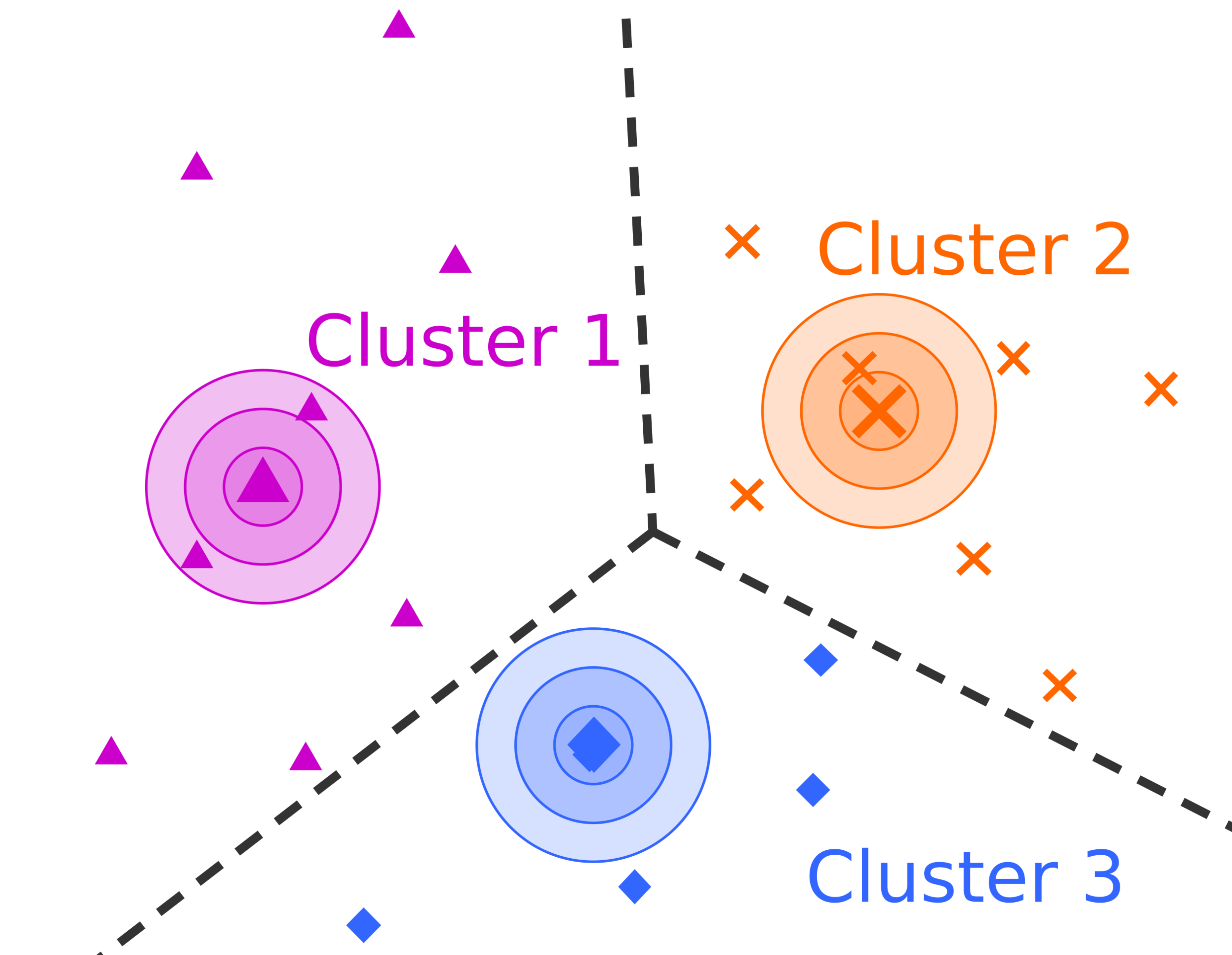}}
	\hfil
	\subfloat[Semi-supervised \texttt{MSSC}]{\includegraphics[width=0.4\textwidth]{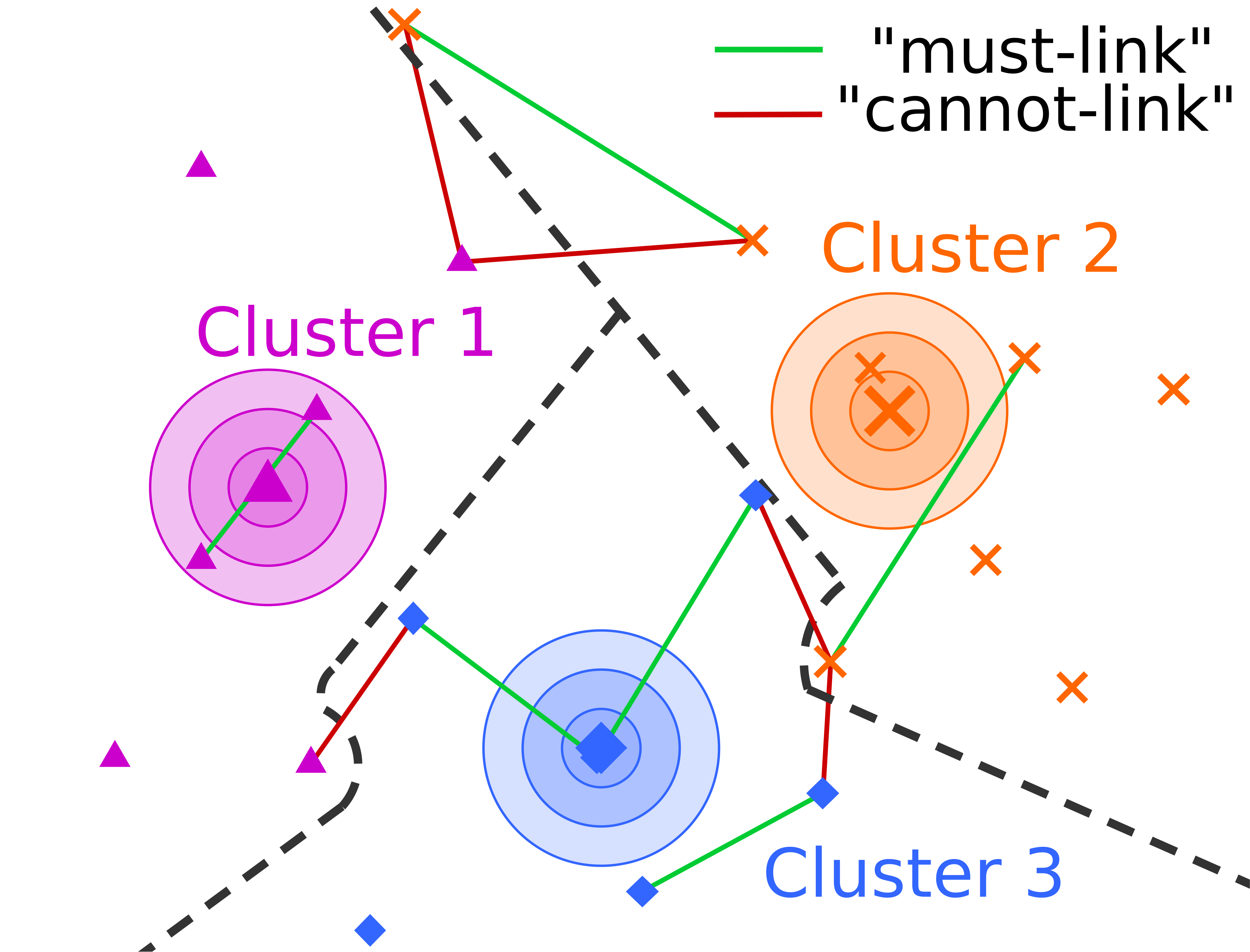}}
	\caption{Two-dimensional examples of assignment step outcome with $N = 20$ points and $K = 3$ cluster centers. The dotted black lines indicate the separations between clusters, while shaded areas just emphasize cluster centers.}
	\label{fig::unc_vs_con}
\end{figure}

It is well known that problem \eqref{eq::mssc-prob} is NP-hard \cite{davidson07}, as it contains the unconstrained \texttt{MSSC} model \cite{aloise_np-hardness_2009} as a special case. Both problems are usually characterized by a large number of local optima which are not global \cite{JAIN2010651} and, thus, they are difficult to solve exactly in terms of both complexity and computational resources. The additional difficulties of the semi-supervised case are mainly derived by the following fact: while given the optimal $x^\star$ of problem \eqref{eq::mssc-prob} we can simply calculate each $k$-th cluster center as $y^\star_k = \frac{1}{\sum_{i = 1}^{N} x^\star_{ik}}\sum_{i = 1}^{N} x^\star_{ik} p_i$ \cite{lloyd82}, the assignment operation, i.e., the mechanism by which we assign each point to a center ($y^\star \rightarrow x^\star$), requires different effort in the unsupervised and semi-supervised scenarios. In the first one, this step is indeed trivial to perform, assigning each point to its closest cluster center; in the semi-supervised context, ``must-link'' and ``cannot-link'' constraints must be taken into account and, then, it is no more true that points are associated to their closest centroid (see Figure \ref{fig::unc_vs_con}).

One of the best known and used technique for \texttt{MSSC} is \texttt{K-MEANS} \cite{hartigan79}. Given an initial position of the cluster centers, \texttt{K-MEANS} iterates the two following steps until cluster centers locations are not updated for two consecutive iterations: (i) each data point is assigned to its closest cluster center; (ii) each center is moved to the arithmetic mean position of the points assigned to its cluster. The first most notable \texttt{K-MEANS} extension for semi-supervised \texttt{MSSC} is \texttt{COP-K-MEANS} \cite{wagstaff2001constrained}, which aims to assign each point to the nearest cluster center such that no constraint is violated; the approach is not equipped with any backtracking procedure and, thus, it may stop failing to return a feasible solution. In order to address this issue, modified versions of \texttt{COP-K-MEANS} have been proposed (e.g., \cite{TAN2010}), with some of the most recent ones \cite{Baumann20, baumann2023pccc} equipped with assignment steps based on binary programming problems to enforce the solution to be feasible at each iteration. Although being quite easy to implement and efficient, \texttt{K-MEANS} usually returns a \textit{local} optimum, which might be very far from the global one. Thus, other paradigms have been explored over the years in order to overcome such limitation: soft-constrained algorithms, i.e., algorithms where penalties for the constraints violations are employed \cite{basu04, Davidson05}, non-smooth optimization approaches \cite{Bagirov23} and exact algorithms \cite{Guns16, PICCIALLI2022105958}. 

In this short manuscript, we propose an extension for semi-supervised \texttt{MSSC} problems of the differential evolution-based memetic framework (\texttt{MDEClust}) proposed in \cite{MANSUETO2021107849}, which must strictly satisfy the given ``must-link'' and ``cannot-link'' constraints. A memetic approach is a population-based meta-heuristic which, at each iteration, maintains and updates a population of candidate solutions through the repeated application of \textit{genetic} operators (e.g., crossover, mutation and selection) in order to return the best generated solution at the end of the execution; what distinguishes it from other evolutionary genetic algorithms is the employment of a local solver for refining the solutions before adding them to the population. As far as we know, our extension, called \texttt{S-MDEClust}, is the first memetic approach in the semi-supervised clustering literature for finding optimal solutions satisfying all the link constraints; the memetic approaches presented so far in the literature \cite{GONZALEZALMAGRO2020104979, GRIBEL2022441} are indeed soft-constrained and, thus, with no guarantee of returning a \textit{feasible} partition.

The main novelties proposed in \texttt{S-MDEClust} are: (i) two novel (exact and greedy) assignment methodologies for the memetic paradigm; (ii) the possible feasibility relaxation carried out by the new greedy assignment algorithm, which turns out to be faster and can drive the algorithm towards the exploration of a wider solution space; (iii) a re-designed mutation operator to deal with the link constraints as well as with the chance of returning a possibly infeasible solution; (iv) the inclusion of the (local) optimization paradigm proposed in \cite{Baumann20} in a memetic scheme.

The remainder of the paper is structured as follows. In Section \ref{sec::preliminaries}, we briefly describe the memetic framework \texttt{MDEClust}. In Section \ref{sec::MDEforSS}, we describe how to extend the latter for semi-supervised \texttt{MSSC}, highlighting the main differences w.r.t.\ the original scheme. In Section \ref{sec::exp_num}, we report the results of thorough computational experiments to highlight the goodness of our proposal. Finally, in Section \ref{sec::conclusions}, we provide some concluding remarks.

\section{Preliminaries}
\label{sec::preliminaries}

In this section, we briefly describe \texttt{MDEClust}, the \textit{Memetic Differential Evolution} (\texttt{MDE}) framework proposed in \cite{MANSUETO2021107849}; the scheme is reported in Algorithm \ref{alg::MDEforSS}.

\begin{algorithm}[h]
	Input: $\mathcal{D}$, $K \in \mathbb{N}^+$, $P \in \mathbb{N}^+$, $N_{MAX} \in \mathbb{N}^+$, $\delta \in \mathbb{R}^+$.\\
	Initialize population $\mathcal{P} = \{S_1,\ldots,S_P\}$ \label{line::initialization}\\
	Let $S^\star \in \mathcal{P}$ be such that $f(S^\star) \le f(S)$, $\forall S \in \mathcal{P}$, $n_{it}^\star = 0$\\
	\While{$n_{it}^\star < N_{MAX}$ $\land$ $\sum_{s=1}^{P}\sum_{\bar{s}>s}^{P}|f(S_s) - f(S_{\bar{s}})| > \delta$ \label{line::stopping-criteria}}{
		\ForAll{$S_s \in \mathcal{P}$}{
			Randomly select $S_1$, $S_2$, $S_3 \in \mathcal{P}$, all different from each other and from $S_s$\\
			Execute crossover with $S_1$, $S_2$, $S_3$ to generate an offspring $O_s$\label{line::crossover}\\
			(Optional) Execute mutation on $O_s$ to get $\tilde{O}_s$\label{line::mutation}\\
			Apply local search to $\tilde{O}_s$ to obtain $O'_s$\label{line::local_search}\\ 
			\If{$f(O'_s) < f(S_s)$ \label{line::if_better_1}}{
				Set $S_s = O'_s$\\
				\If{$f(O'_s) < f(S^\star)$}{
					Set $S^\star = O'_s$, $n_{it}^\star = 0$
				}
				\Else{
					Set $n_{it}^\star = n_{it}^\star + 1$ \label{line::if_better_2}
				}
			}
		}
	}
	\Return $S^\star$
	\caption{\texttt{MDEClust} Framework}
	\label{alg::MDEforSS}
\end{algorithm}

The \texttt{MDEClust} framework initializes and keeps updated throughout the iterations a population of $P$ \texttt{MSSC} solutions $\mathcal{P} = \{S_1, \ldots, S_P\}$. Each solution $S \in \mathcal{P}$ is identified as $S = (\phi^S, \psi^S)$: the \textit{Membership vector} $\phi^S \in \mathbb{N}^N$ is such that, for all $i \in \{1,\ldots, N\}$, $\phi^S_i = k$, with $k \in \{1,\ldots,K\}$, if and only if $x_{ik} = 1$ (see problem \eqref{eq::mssc-prob}); the \textit{Coordinate matrix} $\psi^S = [y_1, y_2,\ldots, y_K]^\top \in \mathbb{R}^{K \times d}$ contains the coordinates of the cluster centers. It is trivial to observe that each of the two items is sufficient to completely characterize a solution, being derivable from each other in $O(NKd)$ in the unsupervised case \cite{GRIBEL2019569, MANSUETO2021107849} similarly to $x^\star$ and $y^\star$ (problem \eqref{eq::mssc-prob}, Section \ref{sec::introduction}). The \texttt{MSSC} objective function value of a solution $S$ is consequently indicated as $f(S) = \sum_{i=1}^{N}\|p_i - \psi^S_{\phi^S_i}\|^2.$

In \texttt{MDEClust}, the initialization of the population is accomplished executing $P$ independent runs of the \texttt{K-MEANS} algorithm \cite{lloyd82} (Line \ref{line::initialization}), in each of which the initial centers are randomly selected among the points of the dataset $\mathcal{D}$ at hand. The population is then updated at each iteration performing the following operations for each solution $S_s \in \mathcal{P}$, with $s \in \{1,\ldots, P\}$: (i) an \textit{offspring} solution $O_s$ is generated by the \textit{crossover} operator (Line \ref{line::crossover}); (ii) the (optional) \textit{mutation} operator is performed, returning a modified offspring $\tilde{O}_s$ (Line \ref{line::mutation}); (iii) \texttt{K-MEANS} is executed as local search procedure to refine $\tilde{O}_s$ (Line \ref{line::local_search}); (iv) the resulting solution $O'_s$ is then compared with $S_s$ in terms of \texttt{MSSC} objective function (see Problem \ref{eq::mssc-prob}) and, if it is better, it replaces $S_s$ in the population (Lines \ref{line::if_better_1}-\ref{line::if_better_2}). The algorithm stops when one of the following stopping criteria is met (Line \ref{line::stopping-criteria}), returning the best \texttt{MSSC} solution in the population: $N_{MAX} > 0$ consecutive iterations have been performed without any best solution improvement; the difference in terms of \texttt{MSSC} objective function among the solutions in the population has fallen below a threshold $\delta > 0$, i.e., $\sum_{s = 1}^{P}\sum_{\bar{s} > s}^{P} | f(S_s) - f(S_{\bar{s}}) | \le \delta$, meaning that the population has collapsed to a single solution. Other stopping criteria can be also possible.

\paragraph{Crossover Operator}

Given a solution $S_s \in \mathcal{P}$, with $s \in \{1,\ldots, P\}$, in \texttt{MDEClust} the \textit{crossover} operator selects from the population three distinct random solutions $S_1, S_2, S_3$, all different from $S_s$, and it creates an offspring solution $O_s = (\phi^{O_s}, \psi^{O_s})$ combining them in the following way: $\psi^{O_s} = \psi^{S_1} + F(\psi^{S_2} - \psi^{S_3})$, with $F \in (0, 2)$.  If $F \approx 0$, a more intensive exploration of $S_1$ neighborhood is carried out; if $F \approx 2$, we likely take the exploration to new regions of problem \eqref{eq::mssc-prob}. The membership vector $\phi^{O_s}$ is then obtained performing an assignment step. 

Since any permutation of cluster centers gives rise to the same clustering solution and no ordering of the cluster centers is a-priori provided, similarly to \cite{GRIBEL2019569} the linear combination is carefully designed based on a bi--partite center matching strategy; the latter allows the re-labeling of cluster centers so that if two solutions, e.g., $S_2$ and $S_3$, are very similar one to the other, their difference $\psi^{S_2} - \psi^{S_3}$ will be close to zero. The exact bi--partite matching strategy employed in this context is the Hungarian algorithm \cite{Kuhn55}, whose complexity is $O(K^3)$.

\paragraph{Mutation Operator}

The \textit{mutation} operator aims to diversify the population, preventing a premature stop of the framework caused by a rapid collapse of the population into a single solution. In \texttt{MDEClust}, the operator mainly consists of a random relocation of one center of the offspring solution $O_s$ generated during the crossover phase, thus generating a modified solution $\tilde{O}_s$ characterized by new features that are not inherited from the parents. More details on the topic can be found in \cite{MANSUETO2021107849}. 

\section{A \texttt{MDE} Framework for Semi-Supervised \texttt{MSSC}}
\label{sec::MDEforSS}

In this section, we propose \texttt{S-MDEClust}, a possible extension of the \texttt{MDEClust} framework for handling semi-supervised \texttt{MSSC} problems. Here, we remark that, to the best of our knowledge, \texttt{S-MDEClust} is the first memetic approach in the semi-supervised clustering literature designed to generate a (hopefully optimal) \textit{feasible} solution.

In the next sections, we describe the main differences w.r.t.\ the original framework: the assignment step, no more trivial as in the unsupervised scenario (see Figure \ref{fig::unc_vs_con}) and the operations highly depending on it, i.e., mutation (Line \ref{line::mutation}) and local search phase (Line \ref{line::local_search}). 

\subsection{Assignment Step}
\label{subsec::assignment_step}

Given a solution $S$ and its coordinate matrix $\psi^S$, we propose two (exact and greedy) assignment steps capable of dealing with link constraints to get the membership vector $\phi^S$.

\paragraph{Exact Assignment Step}

Inspired by \cite{Baumann20}, the exact assignment step is based on solving the optimization problem
\begin{equation}
	\label{eq::asignment-step}
	\begin{aligned}
		\min_{x \in \mathcal{X} \cap \mathcal{L}}\;\sum_{i = 1}^{N}\sum_{k = 1}^{K} x_{ik}\|p_i - \psi^S_k\|^2
	\end{aligned},
\end{equation}
with $\mathcal{X}$ and $\mathcal{L}$ defined as in Equations \eqref{eq::feasible_set_X}-\eqref{eq::feasible_set_L}.

Solving an additional optimization problem each time an assignment step is required surely leads to a greater consumption of computational resources. Problem \eqref{eq::asignment-step} can indeed be seen as a generalization of the minimum cost perfect matching problem with conflict pair constraints \cite{ONCAN2013920}, which is known to be NP--hard. Thus, in order to reduce the computational effort, we also propose a faster greedy assignment methodology, to be used in specific phases of the \texttt{S-MDEClust} framework. 

\paragraph{Greedy Assignment Step}

We report the scheme of the greedy methodology in Algorithm \ref{alg::GAS}.

\begin{algorithm}[!h]
	Input: $\mathcal{D}$, $K \in \mathbb{N}^+$, $\mathcal{M}_\mathcal{L}, \mathcal{C}_\mathcal{L}$, $S=(\phi^S, \psi^S)$.\\
	Let $\mathcal{G} = \left\{G \subseteq \{1,..., N\} \mid \forall i, j \in G \text{ s.t. } i \ne j, (i, j) \in \mathcal{M}_\mathcal{L}\right\}$\label{line::groups}\\
	Let $\mathcal{A} = \emptyset$\\
	\ForAll{$G \in \mathcal{G}$\label{line::listG}}{
		Let $\hat{K} := \{k \in \{1,...,K\} \mid \forall\, (\hat{G}, k) \in \mathcal{A}, i_c \in G, j_c \in \hat{G} \implies (i_c, j_c) \not\in \mathcal{C}_\mathcal{L}\}$\label{line::first_kg}\\
		\textbf{if} $\hat{K} = \emptyset$ \textbf{then} \text{Let } $\hat{K} := \{1,\ldots,K\}$\label{line::nwd}\\
		Let $k_G \in \argmin_{k \in \hat{K}}\,\sum_{i \in G}\|p_i - \psi^S_k\|^2$\label{line::second_kg}\\
		Set $\mathcal{A} = \mathcal{A} \cup \{(G, k_G)\}$\\
		Set $\phi^S_i = k_G$, $\forall i \in G$ \label{line::a}
	}
	\Return $S$
	\caption{Greedy Assignment Step}
	\label{alg::GAS}
\end{algorithm} 

First, all the points are collected in groups based on the ``must-link'' constraints (Line \ref{line::groups}). Then, for each group $G$, we try to find the cluster centers that can be assigned to the points of the group (Line \ref{line::first_kg}): these centers should not be already associated to a group $\hat{G}$ that has a point involved in a ``cannot-link'' constraint with a point of $G$. The methodology then highly depends on the assignment order of the groups, and it could happen that, in presence of ``cannot-link'' constraints, finding a cluster center to associate to a group is not possible. In such cases, we relax the solution feasibility requirement, allowing any cluster center to be assigned to the points of the group (Line \ref{line::nwd}), regardless the possible constraint violations. Given the available cluster centers, we finally assign the group points to one of the nearest (Lines from \ref{line::second_kg} to \ref{line::a}). 

Although Algorithm \ref{alg::GAS} has no \textit{feasibility} guarantee, it is less computationally demanding than the exact assignment step (no binary programming problem is involved). It can then be useful to reduce the computational cost of the internal \texttt{S-MDEClust} operations, i.e., crossover and mutation (Lines \ref{line::crossover}-\ref{line::mutation} of Algorithm \ref{alg::MDEforSS}), where, unlike the local search procedure (Line \ref{line::local_search}), generating a feasible solution is not mandatory. Moreover, starting local searches from non-feasible solutions can be seen as a further mechanism to explore new regions of the feasible set. 

\subsection{Mutation Operator}
\label{subsec::mutation}

The scheme of mutation is reported in Algorithm \ref{alg::mutation_operator}.

\begin{algorithm}[!h]
	Input: $\mathcal{D}$, $K \in \mathbb{N}^+$, $\mathcal{M}_\mathcal{L}, \mathcal{C}_\mathcal{L}$, $O_s = (\phi^{O_s}, \psi^{O_s})$, $\mathcal{U}(\cdot)$ uniform random number generator, $\alpha \in [0, 1]$.\label{line::mut_inputs}\\
	Let $\psi^{\tilde{O}_s} = \psi^{O_s}$, $\bar{c} = \mathcal{U}(1,\ldots, K)$\label{step::center_to_be_removed}\\
	Perform assignment step with $\mathcal{Y} \setminus \{y_{\bar{c}}\}$ to get $\tilde{\phi}^{O_s}$ \label{line::ass_1}\\
	\textbf{if} \textit{assignment step has failed} \textbf{then} Set $\alpha = 0$\label{line::alpha}\\
	Let, $\forall i \in \{1,\ldots,N\}$, $P_i = \frac{1 - \alpha}{N} + \left(\alpha\left\|p_{i} - y_{\tilde{\phi}_{i}^{O_s}}\right\|\right) / \left(\sum_{j = 1}^{N}\left\|p_j - y_{\tilde{\phi}_j^{O_s}}\right\|\right)$ \label{step::prob}\\
	Choose $\bar{\i} \in \{1,\ldots, N\}$ by \textit{roulette wheel} based on $P$ \label{step::roulette}\\ 
	Set $\psi^{\tilde{O}_s}_{\bar{c}} = p_{\bar{\i}}$, Perform assignment step to get $\phi^{\tilde{O}_s}$ \label{line::ass_2}\\
	\Return{$\tilde{O}_s = (\phi^{\tilde{O}_s}, \psi^{\tilde{O}_s})$}
	\caption{Mutation Operator}
	\label{alg::mutation_operator}
\end{algorithm} 

As mentioned in Section \ref{sec::preliminaries}, the operator consists in a random-based relocation of a center of the offspring solution $O_s$ generated by the crossover operator: (i) a center randomly selected with uniform probability is removed (Line \ref{step::center_to_be_removed}); (ii) each dataset point is re-assigned to its closest center among the $K - 1$ remaining ones, returning a new temporary membership vector $\tilde{\phi}^{O_s}$ (Line \ref{line::ass_1}); (iii) the missing center is now restored as a point $p_{\bar{\i}}$, with $\bar{\i} \in \{1,\ldots, N\}$, randomly chosen by \textit{roulette wheel} (Line \ref{step::roulette}), a randomized operation to select potentially useful points based on the probability $P$ (Line \ref{step::prob}). If $\alpha = 0$, every point has the same chances to be selected, thus potentially leading to the exploration of new solutions; as $\alpha$ gets larger, the selection is biased towards points that are more distant from their closest centers.

In Line \ref{line::ass_2}, the assignment step is carried out smoothly, either in the exact or the greedy way (Section \ref{subsec::assignment_step}). On the other hand, in Line \ref{line::ass_1} we have $K-1$ centers to which the points can be assigned; thus, minor modifications in the exact and the greedy methodology are required. In the exact assignment step, problem \eqref{eq::asignment-step} is solved with the second set of constraints defining $\mathcal{X}$ changed into $\sum_{i = 1}^{N} x_{ik} \ge 1, \forall k \in \{1,\ldots, K\} \setminus \{\bar{c}\},\ \sum_{i = 1}^{N} x_{i\bar{c}} = 0$. If the semi-supervised \texttt{MSSC} problem presents some ``cannot-link'' constraints, the \textit{modified} problem \eqref{eq::asignment-step} may not have a feasible solution and, as a consequence, the exact assignment step may fail; being the point-center distances not available in this case, we set $\alpha = 0$ in Line \ref{line::alpha}, thus giving to each dataset point the same chances to be selected by \textit{roulette wheel}. As for the greedy methodology, in Line \ref{line::listG} of Algorithm \ref{alg::GAS} we only take into account the sets $G \in \mathcal{G}$ whose points are assigned to the removed center $y_{\bar{c}}$ (i.e., $\forall i \in G,\ \phi_i^{O_s} = \bar{c}$); moreover, we add the constraint $k_G \ne \bar{c}$ in Line \ref{line::second_kg} so that the new center chosen for a group cannot be the removed one. Unlike the exact assignment step, the greedy methodology always returns a solution, although it might be infeasible.

\subsection{Semi-supervised Local Search Procedure}
\label{subsec::ssk}

In the \texttt{S-MDEClust} framework, we employ as local search procedure the semi-supervised \texttt{K-MEANS} variant proposed in \cite{Baumann20}, which is a standard \texttt{K-MEANS} algorithm (see Section \ref{sec::introduction}) where the assignment step is performed by \textit{exactly} solving problem \eqref{eq::asignment-step}. This local search is used both for the population initialization and the solutions refinement (Lines \ref{line::initialization}-\ref{line::local_search} of Algorithm \ref{alg::MDEforSS}), always returning a feasible solution and, thus, guaranteeing the feasibility of all the solutions inserted in the population throughout the \texttt{S-MDEClust} iterations.

\section{Numerical Experiments}
\label{sec::exp_num}

\begin{table*}[!h]
	\centering
	\renewcommand{\arraystretch}{1.5}%
	\tiny
	\setlength\tabcolsep{2.5pt}
	\begin{minipage}{0.925\textwidth}
		\begin{minipage}{0.49\textwidth}
			\centering
			\begin{tabular}{|c||ccc||ccc||c|}
				\hline
				\multirow{2}{*}{\textit{Dataset}}&\multirow{2}{*}{$N$}&\multirow{2}{*}{$d$}&\multirow{2}{*}{$K$}&\multicolumn{3}{c||}{\textit{\% Viol. (Mean)}}&\multirow{2}{*}{$(|\mathcal{M}_\mathcal{L}|, |\mathcal{C}_\mathcal{L}|)$}\\
				&&&&$\mathcal{M}_\mathcal{L}$&$\mathcal{M}_\mathcal{L}$ \& $\mathcal{C}_\mathcal{L}$&$\mathcal{C}_\mathcal{L}$&\\
				\hline
				\hline
				Iris&150&4&3&16.1&14.9&9.7&\multirow{5}{60pt}{\centering (50, 0), (100, 0), (25, 25), (50, 50), (0, 50), (0, 100)}\\
				\cline{1-7}
				Wine&178&13&3&41.3&32.8&20.8&\\
				\cline{1-7}
				Connectionist&208&60&2&51.6&47.2&51.6&\\
				\cline{1-7}
				Seeds&210&7&3&18.4&16.0&10.7&\\
				\cline{1-7}
				Glass&214&9&6&57.7&40.9&24.0&\\
				\hline
				Accent&329&12&6&76.7&46.3&18.7&\multirow{2}{70pt}{\centering (100, 0), (150, 0), (50, 50), (75, 75), (0, 100), (0, 150)}\\
				\cline{1-7}
				Ecoli&336&7&8&57.9&31.9&3.4&\\
				\hline
				ECG5000&500&140&5&50.0&27.7&5.2&\multirow{2}{70pt}{\centering (150, 0), (250, 0), (75, 75), (125, 125), (0, 150), (0, 250)}\\
				\cline{1-7}
				Computers&500&720&2&46.6&49.8&52.3&\\
				\hline
				\multirow{2}{*}{Gene}&\multirow{2}{*}{801}&\multirow{2}{*}{20531}&\multirow{2}{*}{5}&\multirow{2}{*}{1.3}&\multirow{2}{*}{0.9}&\multirow{2}{*}{0.5}&\multirow{2}{75pt}{\centering (200, 0), (400, 0), (100, 100), (200, 200), (0, 200), (0, 400)}\\
				&&&&&&&\\
				\hline
			\end{tabular}
			\begin{tabular}{|c|c|c|c||c|c|c|c|}%
				\hline%
				\multirow{2}{*}{\textit{Dataset}}&\multirow{2}{*}{$N$}&\multirow{2}{*}{$d$}&\multirow{2}{*}{$K$}&\multicolumn{4}{c|}{\textit{\%Viol. (Mean)}}\\%
				\cline{5%
					-%
					8}%
				&&&&100&200&500&1000\\%
				\hline%
				\hline%
				Banana&5300&2&2&43&55&50&50.1\\%
				Letter&20000&16&26&7.7&4.8&7.2&6.1\\%
				Shuttle&57999&9&7&26.7&32.3&37.3&37\\%
				Mnist&70000&784&10&12.7&11.3&11.5&14.8\\%
				Cifar10&60000&3072&10&17&18.5&15.9&17.6\\%
				\hline%
			\end{tabular}%
		\end{minipage}
		\hfill
		\begin{minipage}{0.49\textwidth}
			\centering
			\begin{tabular}{|c|c||c|c|c|c|c|c|c|}%
				\hline%
				\multirow{2}{*}{$N$}&\multirow{2}{*}{$K$}&\multicolumn{7}{c|}{\textit{\%Viol. (Mean)}}\\%
				\cline{3%
					-%
					9}%
				&&100&200&500&1000&2000&5000&10000\\%
				\hline%
				\hline%
				500&2&58&45.5&48.8&49.2&50.7&50.6&49.8\\%
				500&20&16.3&9.2&11.4&9.6&10.2&9.6&9.3\\%
				500&100&1.7&2.7&1.1&2&2&2&2\\%
				\hline%
				1000&2&41.3&45.2&45.7&40.5&43.6&41.9&42.3\\%
				1000&20&9.3&10.3&7.5&9.3&10.9&9.5&10.5\\%
				1000&100&3.7&0.8&2.5&2.5&1.8&2&2.1\\%
				\hline%
				2000&2&43.3&44.8&44.9&44.9&43.3&47&44.1\\%
				2000&20&4&6.8&13.2&9.6&9.9&10.3&10.1\\%
				2000&100&2.7&4.2&1.7&1.9&1.6&2&2\\%
				\hline%
				5000&2&45&46.5&42.2&45.4&46.1&44.8&44.7\\%
				5000&20&12.7&9.5&10.5&9.6&9.6&9.4&10.6\\%
				5000&100&0.3&3.8&1.3&1.7&2&2&2.2\\%
				\hline%
				10000&2&43&44.8&46&46&46.1&45.8&45.9\\%
				10000&20&12&12&10.3&9.9&10.4&10.2&10.2\\%
				10000&100&2.7&2.5&1.4&2.1&2&2.3&2.4\\%
				\hline%
				20000&2&39.7&43&43.5&42.7&44.1&43.9&45\\%
				20000&20&9&11&10.9&8.5&10.7&10.3&10.4\\%
				20000&100&0&2.2&2.3&1.6&2.4&2.5&2.1\\%
				\hline%
			\end{tabular}%
		\end{minipage}
	\end{minipage}
	\caption{$1^{\text{st}}$ set: well-known datasets with constraint configurations; for each dataset and constraint type, the average percentage of violated pairwise constraints in the global minimum of the unconstrained \texttt{MSSC} problem is reported from \cite{PICCIALLI2022105958}. $2^{\text{nd}}$ set: large-size datasets. $3^{\text{rd}}$ set: two-dimensional synthetic instances generated by mixtures of spherical Gaussian distributions. For each dataset of the $2^{\text{nd}}$ and $3^{\text{rd}}$ sets and for each number of constraints, the average percentage of violated pairwise constraints obtained by 3 runs of \texttt{K-MEANS} solving the unconstrained \texttt{MSSC} problem is provided.}
	\label{tab::datasets}
\end{table*}

In this section, we report the results of thorough computational experiments, where we compared our framework \texttt{S-MDEClust} with some state-of-the-art methodologies from the semi-supervised clustering literature. The \texttt{S-MDEClust} code was written in \texttt{Python3}. All the tests were run on a computer with the following characteristics: Intel Xeon Processor E5-2430 v2, 6 cores, 12 threads, 2.50 GHz, 32 GB RAM. In order to solve instances of problem \eqref{eq::asignment-step}, the Gurobi Optimizer \cite{gurobi} (Version 11) was employed.

\texttt{S-MDEClust} was tested in four different variants: \texttt{S-MDE} (no mutation and exact assignment), \texttt{SM-MDE} (mutation and exact assignment), \texttt{SG-MDE} (no mutation and greedy assignment) and \texttt{SMG-MDE} (mutation and greedy assignment). We remind that in the local search procedure (Section \ref{subsec::ssk}) the exact assignment step is always employed. We set the parameters values according to the \texttt{MDEClust} original paper \cite{MANSUETO2021107849} and some preliminary experiments, not reported here for the sake of brevity: $P=20$; $N_{MAX} = 500$; $\delta = 10^{-4}$; $F \in [0.5, 0.8]$, randomly chosen at each crossover execution; $\alpha=0.5$ for the mutation operator.

The datasets used for experimentation are listed in Table \ref{tab::datasets}. The $1^{\text{st}}$ set is composed by well-known datasets of different sizes and/or number of features, which can be downloaded from the UCI (\href{https://archive.ics.uci.edu}{archive.ics.uci.edu}) and UCR (\href{https://www.cs.ucr.edu/~eamonn/time_series_data_2018/}{cs.ucr.edu}) repositories. For each dataset, we considered the 30 constraint configurations\footnote{The constraint configurations for the $1^{\text{st}}$ dataset set (Table \ref{tab::datasets}) can be downloaded at \href{https://github.com/antoniosudoso/pc-sos-sdp}{github.com/antoniosudoso/pc-sos-sdp.}} generated in \cite{PICCIALLI2022105958} according to the true class partitioning. The $2^{\text{nd}}$ and $3^{\text{rd}}$ sets contain large-size\footnote{The large-size datasets can be found in \href{https://github.com/phil85/PCCC-Data}{github.com/phil85/PCCC-Data.}} and synthetic datasets; for each of these datasets, constraints were created randomly selecting pairs of samples and, based on the true class partitioning, defining ``must-link'' or ``cannot-link'' constraints. The synthetic datasets were generated by mixtures of spherical Gaussian distributions, with the standard deviation of the resulting clusters set to 10 so as they did not result well-separated. 

The metrics for the comparisons are the following: \texttt{MSSC} objective function (\textit{MSSC O.F.}), elapsed time (\textit{Time}), number of calls to the local search (\textit{N° L.S.}) and total number of local search iterations (\textit{N° L.S.\ Iters}), which is more suited than \textit{N° L.S.} when comparing structurally different algorithms. In order to compare an \texttt{S-MDEClust} variant and a competitor w.r.t.\ a metric, we also made use of the relative percentage difference defined as $\Delta = \frac{100(m - m_R)}{m_R}$, where $m$ and $m_R$ is the metric values achieved by the variant and the competitor, respectively; a negative $\Delta$ means that the variant outperformed the competitor in terms of that metric. Finally, in order to summarize some numerical results, we made use of the performance profiles \cite{dolan2002benchmarking}, a graphical tool showing the probability that a metric value achieved by a method in a problem is within a factor $\tau \in \mathbb{R}$ of the best value obtained by any of the algorithms in that problem.

\subsection{Preliminary Assessment of the \texttt{S-MDEClust} Variants}
\label{subsec::prel}

\begin{figure*}[!h]
	\centering
	\subfloat{\includegraphics[width=0.25\textwidth]{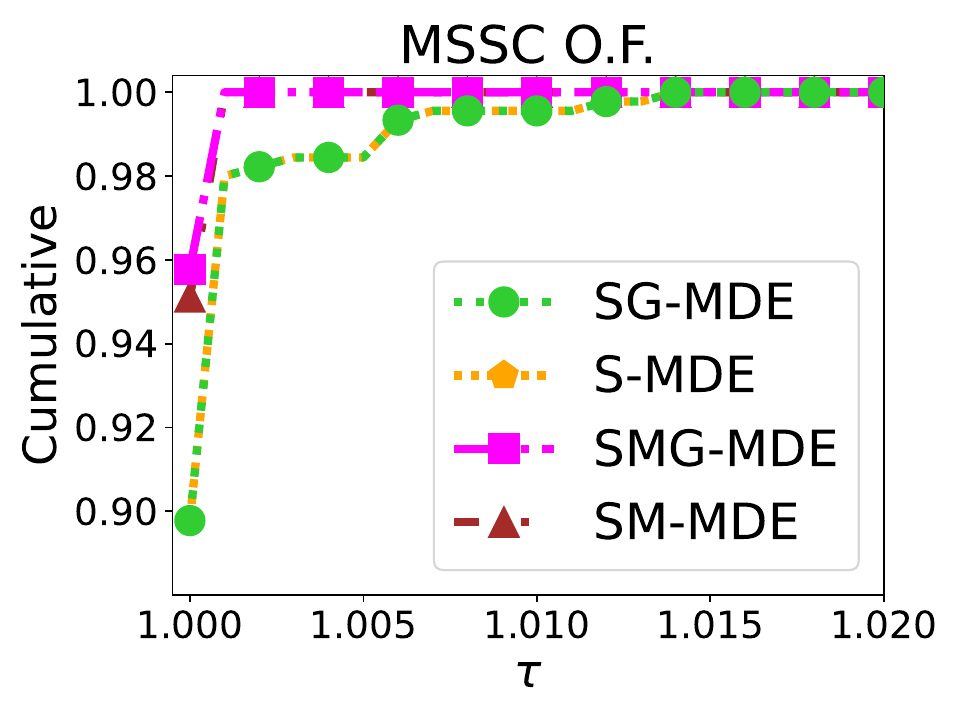}}
	\hfil
	\subfloat{\includegraphics[width=0.25\textwidth]{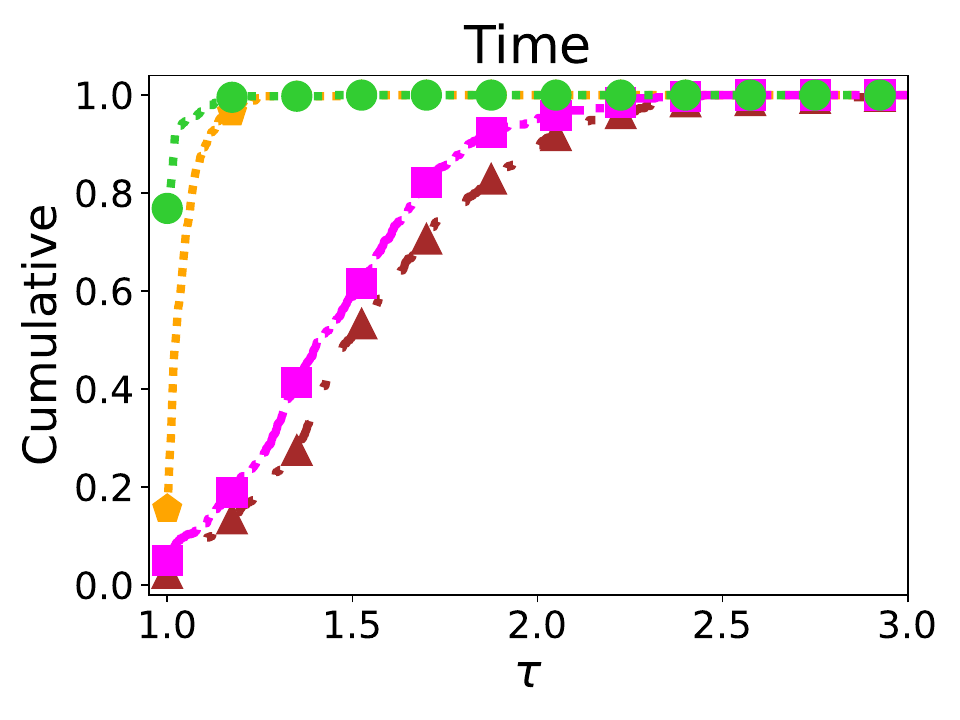}}
	\hfil
	\subfloat{\includegraphics[width=0.25\textwidth]{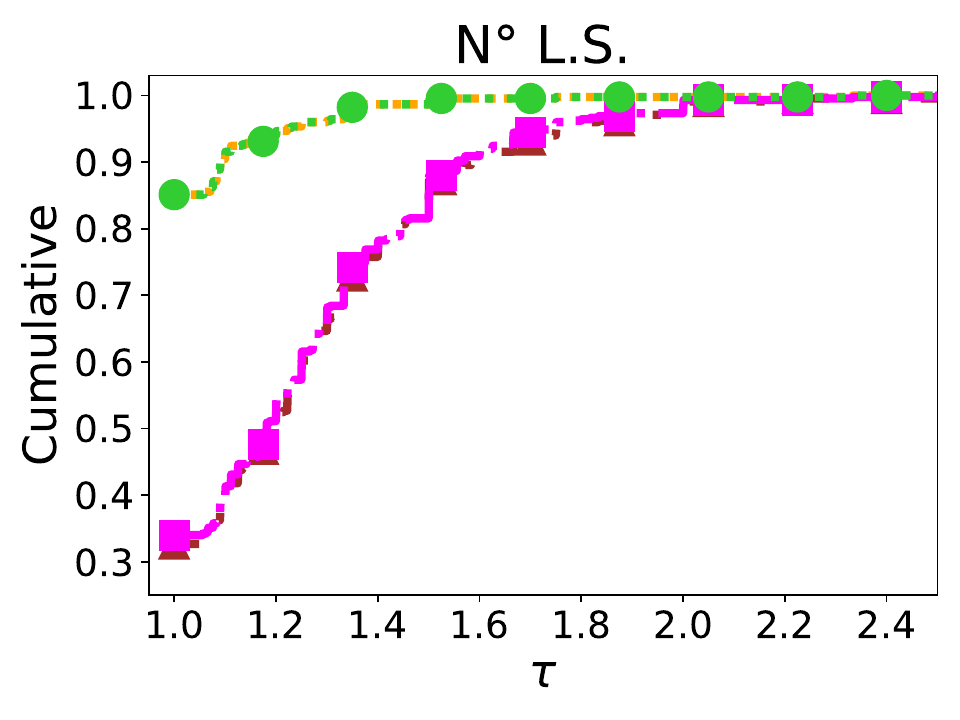}}
	\hfil
	\subfloat{\includegraphics[width=0.25\textwidth]{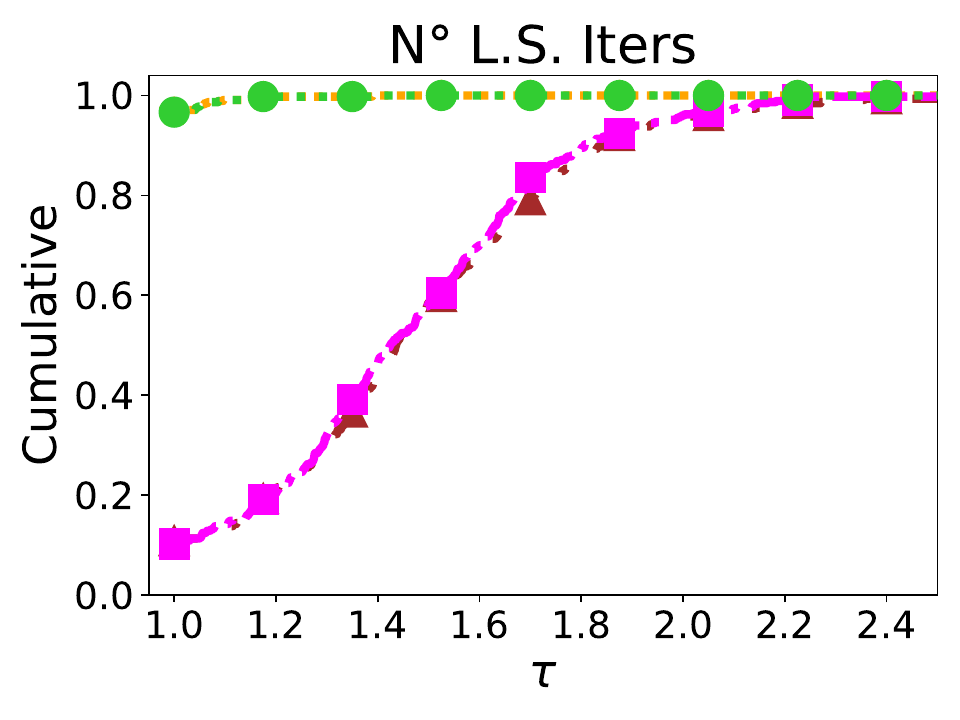}}
	\caption{Performance profiles for \texttt{SG-MDE}, \texttt{S-MDE}, \texttt{SMG-MDE} and \texttt{SM-MDE} on the datasets Iris, Accent and ECG5000 (see Table \ref{tab::datasets}). Note that the intervals of the axes were set for a better visualization of the numerical results.}
	\label{fig::preliminary}
\end{figure*}

We first compared the \texttt{S-MDEClust} variants on a subset of datasets (Iris, Accent, ECG5000), in order to observe the influence of the mutation operator presence and/or the use of the greedy assignment step; in Figure \ref{fig::preliminary}, we report the performance profiles obtained by 5 runs of the variants with 5 different seeds for the pseudo-random number generator. All the variants performed well in terms of \textit{MSSC O.F.}: the best results were obtained only employing the mutation in the framework, at the expense, however, of a greater number of local search executions and, as a consequence, of a greater elapsed time. The greedy assignment step allowed to save computational time as can be observed by the performance profiles in terms of \textit{Time} for \texttt{SG-MDE} and \texttt{S-MDE}, two variants that executed the same number of calls to the local search overall (a similar result can be also observed for \texttt{SMG-MDE} and \texttt{SM-MDE}). Taking into account these results, we decided to consider only \texttt{SG-MDE} for the next comparisons.

\subsection{Comparison with well-known Heuristic Approaches}
\label{subsec::comparisons-heuristics}

Since, as far as we know, no other memetic approach designed to return \textit{strictly feasible} solutions has been proposed in the literature yet, we tested the goodness of \texttt{S-MDEClust} comparing it with well-known \textit{local} heuristics for semi-supervised clustering run in a \textit{multi-start} fashion: \texttt{COP-K-MEANS} \cite{wagstaff2001constrained}, introduced in Section \ref{sec::introduction} and referred as \texttt{COP-KM}; \texttt{BLP-KM} \cite{Baumann20}, the Baumann approach also employed as local search in our framework (Section \ref{subsec::ssk}). The parameters settings for the two competitors were the ones recommended in \cite{wagstaff2001constrained, Baumann20}. \texttt{SG-MDE}, \texttt{COP-KM} and \texttt{BLP-KM} were executed in order to process 100 solutions: each execution of \texttt{COP-KM} and \texttt{BLP-KM} consisted in 100 runs with 100 initial random solutions, while for \texttt{SG-MDE} we considered the population size $P=10$ and a maximum number of iterations equal to 10. For all the algorithms, the maximum number of iterations for the local search was set to 25. All the tests were performed 3 times with different seeds for the pseudo-random number generator; the shown metrics results are the averages of the values obtained in the 3 executions.

The results of the experiments on the synthetic datasets are summarized in Table \ref{tab::general-synthetic}, where we report for each dataset: the percentage (\textit{\%S}) of constraint configurations where the solvers reached the best solution quality; the percentage of times \texttt{SG-MDE} performed the minimum number of local search iterations; the rate \textit{\%}$N_\text{unfeas}$ of times \texttt{COP-KM} did not manage to find any feasible solution. \texttt{SG-MDE} was by far the most effective and efficient approach in the tested benchmark, with the performance difference w.r.t.\ \textit{MSSC O.F.} being clearer as the values for $N$ and $K$ grows. As for \textit{N° L.S. Iters}, the lack of columns for \texttt{BLP-KM} and \texttt{COP-KM} is justified since they did not manage to perform less local search iterations than our approach, which successfully exploited the population convergence based stopping condition to save computational resources without spoiling the solution quality. 

\begin{table}[h]
	\tiny
	\centering
	\hskip-0.2cm
	\renewcommand{\arraystretch}{1.3}%
	\begin{tabular}{|c|c||ccc||c||c|}%
		\hline%
		\multirow{3}{*}{$N$}&\multirow{3}{*}{$K$}&\multicolumn{4}{c||}{\textit{\%S}}&\multirow{2}{*}{\textit{\%}$N_\text{unfeas}$}\\%
		\cline{3%
			-%
			6}%
		&&\multicolumn{3}{c||}{\textit{MSSC O.F.}}&\textit{N° L.S. Iters}&\\%
		\cline{3%
			-%
			7}%
		&&\texttt{SG-MDE}&\texttt{BLP-KM}&\texttt{COP-KM}&\texttt{SG-MDE}&\texttt{COP-KM}\\%
		\hline%
		\hline%
		500&2&\textbf{100}&\textbf{100}&57.14&\textbf{100}&42.86\\%
		500&20&\textbf{100}&14.29&0&\textbf{100}&28.57\\%
		500&100&\textbf{100}&0&0&\textbf{100}&0\\%
		\hline%
		1000&2&\textbf{100}&\textbf{100}&42.86&\textbf{100}&42.86\\%
		1000&20&\textbf{100}&0&0&\textbf{100}&14.29\\%
		1000&100&\textbf{100}&0&0&\textbf{100}&0\\%
		\hline%
		2000&2&\textbf{100}&85.71&28.57&\textbf{100}&57.14\\%
		2000&20&\textbf{100}&0&0&\textbf{100}&0\\%
		2000&100&\textbf{100}&0&0&\textbf{100}&0\\%
		\hline%
		5000&2&\textbf{100}&57.14&0&\textbf{100}&57.14\\%
		5000&20&\textbf{100}&0&0&\textbf{100}&0\\%
		5000&100&\textbf{100}&0&0&\textbf{100}&0\\%
		\hline%
		10000&2&\textbf{100}&71.43&0&\textbf{100}&71.43\\%
		10000&20&\textbf{100}&0&0&\textbf{100}&0\\%
		10000&100&\textbf{100}&0&0&\textbf{100}&0\\%
		\hline%
		20000&2&\textbf{100}&42.86&0&\textbf{100}&71.43\\%
		20000&20&\textbf{100}&0&0&\textbf{100}&0\\%
		20000&100&\textbf{100}&0&0&\textbf{100}&0\\%
		\hline%
	\end{tabular}%
	\caption{Numerical results obtained by \texttt{SG-MDE}, \texttt{BLP-KM} and \texttt{COP-KM} on the synthetic datasets of Table \ref{tab::datasets}.}
	\label{tab::general-synthetic}
\end{table}

\begin{figure}[!h]
	\centering
	\subfloat{\includegraphics[width=0.25\textwidth]{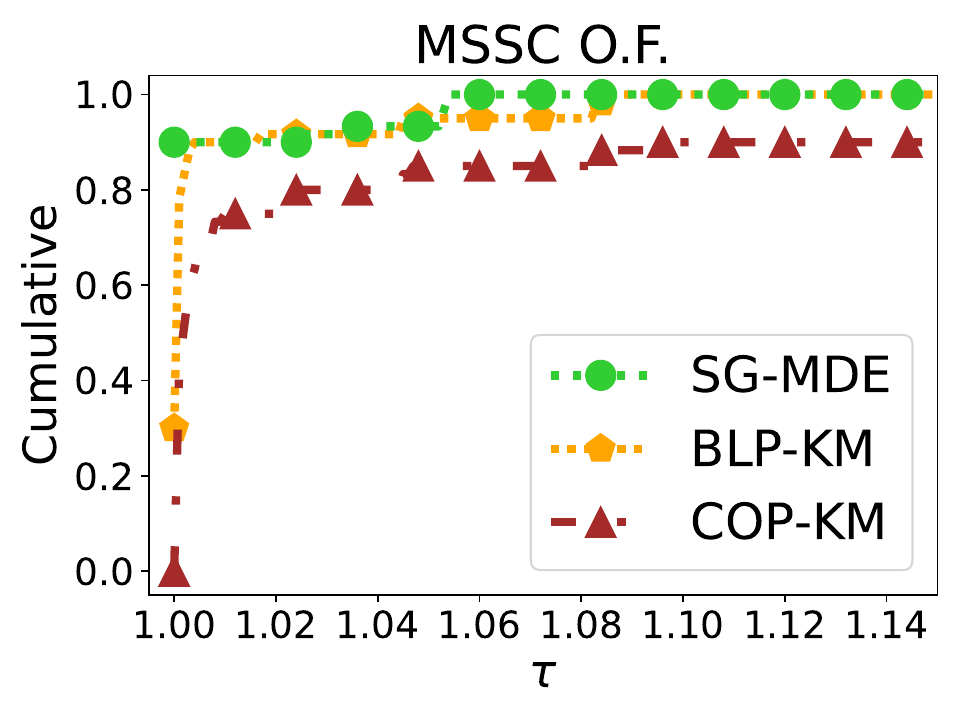}}
	\hfil
	\subfloat{\includegraphics[width=0.25\textwidth]{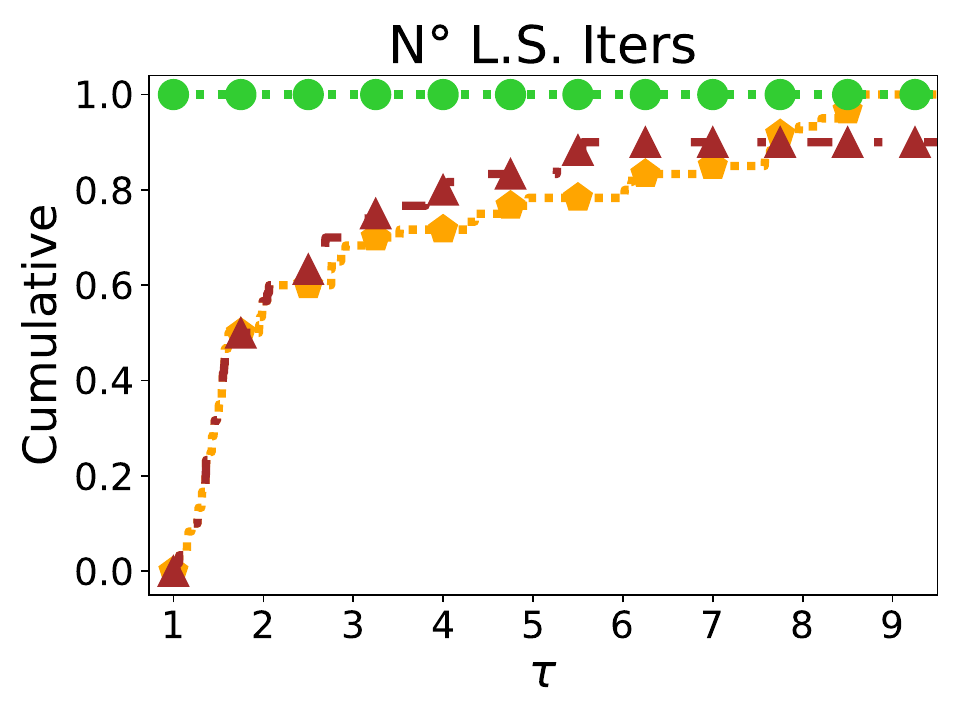}}
	\caption{Performance profiles for \texttt{SG-MDE}, \texttt{BLP-KM} and \texttt{COP-KM} on the large-size datasets ($2^{\text{nd}}$ set of Table \ref{tab::datasets}). Note that the intervals of the x-axes were set for a better visualization of the numerical results.}
	\label{fig::pp-heuristics}
\end{figure}

The performance profiles (Figure \ref{fig::pp-heuristics}) obtained on the large-size datasets of Table \ref{tab::datasets} confirm the superiority of \texttt{SG-MDE} w.r.t. the competitors in terms of both metrics.

\subsection{Comparison with Global State-of-the-art Methodologies}
\label{subsec::comparison_global}

We finally decided to test our framework with some recent, albeit of different types, state-of-the-art global optimization approaches: the \textit{exact} branch-and-cut algorithm proposed in \cite{PICCIALLI2022105958}, called \texttt{PC-SOS-SDP}; the soft-constrained memetic approach presented in \cite{GRIBEL2022441}, which we call \texttt{S-HG-MEANS}. We set the parameters for the 2 competitors according to \cite{PICCIALLI2022105958,GRIBEL2022441}. Since in our setting the constraints were considered accurate, we set the \texttt{S-HG-MEANS} parameter indicating the accuracy of the provided constraints $p=1$. As they involve randomized operations, \texttt{SG-MDE} and \texttt{S-HG-MEANS} were run 5 times, and the results were again obtained by calculating the averages of the values achieved in the 5 tests. In Table \ref{tab::SM-MDEvsPC}, we compare the performance of \texttt{SG-MDE} and \texttt{PC-SOS-SDP} on the $1^{\text{st}}$ dataset set of Table \ref{tab::datasets}. For each dataset we report: \textit{\%S}, the percentage of constraint configurations where \texttt{SG-MDE} obtained a better or equal performance w.r.t.\ \texttt{PC-SOS-SDP}; $\Delta_\text{Mi}$, $\Delta_\text{Me}$ and $\Delta_\text{Ma}$, which are, respectively, the minimum, median and maximum of all the relative percentage differences obtained on the tested constraint configurations; for the \textit{Time} metric, we only considered the instances where our methodology was better in terms of \textit{MSSC O.F.}. Since \texttt{PC-SOS-SDP} spends additional time to prove the solution optimality, in order to have a fairer comparison we decided to consider as \textit{Time} value for the exact approach the elapsed time to find the best solution. As for \textit{MSSC O.F.}, \texttt{SG-MDE} proved to be competitive with \texttt{PC-SOS-SDP}, reaching most of the times a similar (and sometimes even better) performance. The results w.r.t. \textit{Time} are justified by the different nature of the algorithms: the exact approach is composed by more time-consuming operations than our memetic methodology, which, as a consequence, achieved the best performance overall. The \textit{Time} metric is thus inconsistent alone, but the analysis of it in conjunction with \textit{MSSC O.F.} leads to an interesting, not obvious a-priori, conclusion: \texttt{SG-MDE}, with only fractions of the CPU times needed by \texttt{PC-SOS-SDP}, managed to reach an effectiveness level which is highly competitive with the one of an exact algorithm. In Table \ref{tab::SM-MDEvsSSC}, we report a similar comparison with \texttt{S-HG-MEANS} on the dataset instances where no constraint is violated by \texttt{S-HG-MEANS}. Although our approach was better than the competitor in terms of \textit{N° L.S.}, it took more time than \texttt{S-HG-MEANS} in most of the instances; these results are a direct consequence of the fact that \texttt{S-HG-MEANS} employs a less computationally demanding local search, which, however, does not have any feasibility guarantee. Anyway, such efficiency did not allow \texttt{S-HG-MEANS} to reach feasibility in most of the scenarios; moreover, in the instances where it did, it was outperformed (with remarkable results in some cases) by \texttt{SG-MDE}. 

\begin{table*}[!h]
	\centering
	\renewcommand{\arraystretch}{1.5}%
	\tiny
	\setlength\tabcolsep{2.5pt}
	\begin{tabular}{|c||cccc||cccc|}%
		\hline%
		\multirow{2}{*}{\textit{Dataset}}&\multicolumn{4}{c||}{\textit{MSSC O.F.}}&\multicolumn{4}{c|}{\textit{Time}}\\%
		\cline{2%
			-%
			9}%
		&\textit{\%S}&$\Delta_\text{Mi}$&$\Delta_\text{Me}$&$\Delta_\text{Ma}$&\textit{\%S}&$\Delta_\text{Mi}$&$\Delta_\text{Me}$&$\Delta_\text{Ma}$\\%
		\hline%
		\hline%
		Iris&\textbf{100.0}&\textbf{0.0}&\textbf{0.0}&\textbf{0.0}&\textbf{93.33}&\textbf{{-}95.02}&\textbf{{-}70.09}&35.73\\%
		Wine&\textbf{100.0}&\textbf{0.0}&\textbf{0.0}&\textbf{0.0}&\textbf{96.67}&\textbf{{-}98.2}&\textbf{{-}94.01}&59.79\\%
		Connectionist&\textbf{93.33}&\textbf{{-}0.001}&\textbf{0.0}&0.001&\textbf{93.33}&\textbf{{-}98.2}&\textbf{{-}94.58}&\textbf{{-}83.2}\\%
		Seeds&\textbf{100.0}&\textbf{0.0}&\textbf{0.0}&\textbf{0.0}&\textbf{100.0}&\textbf{{-}97.54}&\textbf{{-}89.46}&\textbf{{-}38.69}\\%
		Glass&23.33&\textbf{0.0}&0.014&0.425&23.33&\textbf{{-}93.01}&\textbf{{-}88.64}&\textbf{{-}82.49}\\%
		\hline%
		Accent&46.67&\textbf{{-}0.004}&0.0&0.409&46.67&\textbf{{-}97.3}&\textbf{{-}91.4}&\textbf{{-}81.8}\\%
		Ecoli&36.67&\textbf{0.0}&0.003&0.479&36.67&\textbf{{-}94.75}&\textbf{{-}88.0}&\textbf{{-}54.15}\\%
		ECG5000&46.67&\textbf{0.0}&0.0&0.03&46.67&\textbf{{-}99.37}&\textbf{{-}93.29}&\textbf{{-}72.48}\\%
		Computers&\textbf{76.67}&\textbf{0.0}&\textbf{0.0}&0.001&\textbf{76.67}&\textbf{{-}99.27}&\textbf{{-}98.17}&\textbf{{-}88.56}\\%
		Gene&\textbf{100.0}&\textbf{0.0}&\textbf{0.0}&\textbf{0.0}&\textbf{100.0}&\textbf{{-}92.12}&\textbf{{-}78.99}&\textbf{{-}57.45}\\%
		\hline%
	\end{tabular}%
	\caption{Numerical results obtained by the \texttt{SG-MDE} algorithm w.r.t. \texttt{PC-SOS-SDP} on the $1^{\text{st}}$ dataset set of Table \ref{tab::datasets}.}
	\label{tab::SM-MDEvsPC}
\end{table*}

\begin{table*}[!h]
	\centering
	\renewcommand{\arraystretch}{1.5}%
	\tiny
	\setlength\tabcolsep{2.5pt}
	\begin{tabular}{|c||cccc||cccc||cccc|}%
		\hline%
		\multirow{2}{*}{\textit{Dataset}}&\multicolumn{4}{c||}{\textit{MSSC O.F.}}&\multicolumn{4}{c||}{\textit{Time}}&\multicolumn{4}{c|}{\textit{N° L.S.}}\\%
		\cline{2%
			-%
			13}%
		&\textit{\%S}&$\Delta_\text{Mi}$&$\Delta_\text{Me}$&$\Delta_\text{Ma}$&\textit{\%S}&$\Delta_\text{Mi}$&$\Delta_\text{Me}$&$\Delta_\text{Ma}$&\textit{\%S}&$\Delta_\text{Mi}$&$\Delta_\text{Me}$&$\Delta_\text{Ma}$\\%
		\hline%
		\hline%
		Iris (30)&\textbf{100.0}&\textbf{{-}1.28}&\textbf{{-}0.001}&\textbf{0.0}&40.0&\textbf{{-}52.83}&18.59&160.91&\textbf{86.67}&\textbf{{-}94.45}&\textbf{{-}79.08}&33.33\\%
		Seeds (30)&\textbf{100.0}&\textbf{{-}1.134}&\textbf{{-}0.044}&\textbf{0.0}&26.67&\textbf{{-}60.81}&24.69&211.93&\textbf{90.0}&\textbf{{-}96.08}&\textbf{{-}67.31}&40.0\\%
		Glass (3)&\textbf{100.0}&\textbf{{-}29.552}&\textbf{{-}28.591}&\textbf{{-}26.053}&0.0&3.11&7.23&94.26&\textbf{100.0}&\textbf{{-}75.16}&\textbf{{-}74.51}&\textbf{{-}64.71}\\%
		\hline%
		Accent (10)&\textbf{100.0}&\textbf{{-}2.72}&\textbf{{-}1.73}&\textbf{{-}0.985}&0.0&77.43&143.9&191.1&\textbf{100.0}&\textbf{{-}66.27}&\textbf{{-}52.78}&\textbf{{-}43.8}\\%
		Ecoli (16)&\textbf{100.0}&\textbf{{-}22.062}&\textbf{{-}16.94}&\textbf{{-}4.068}&0.0&69.79&119.63&334.3&\textbf{100.0}&\textbf{{-}58.43}&\textbf{{-}46.05}&\textbf{{-}1.75}\\%
		\hline%
	\end{tabular}%
	\caption{Numerical results obtained by the \texttt{SG-MDE} algorithm w.r.t. \texttt{S-HG-MEANS} on the dataset instances of the $1^{\text{st}}$ set (Table \ref{tab::datasets}) where no constraint is violated by \texttt{S-HG-MEANS}; the numbers of such instances are indicated in the parentheses.}
	\label{tab::SM-MDEvsSSC}
\end{table*}

\section{Conclusions}
\label{sec::conclusions}

In this paper, we proposed \texttt{S-MDEClust}, an extension of the differential evolution based memetic framework proposed in \cite{MANSUETO2021107849} for semi-supervised \texttt{MSSC} problems, where instance-level constraints, each of them indicating the membership of two points to the same or to a different cluster, are given. As far as we know, \texttt{S-MDEClust} represents the first memetic algorithm in the semi-supervised clustering literature designed to return a (hopefully) optimal \textit{feasible} solution. The main modifications w.r.t. the original involves specific operations, such as the assignment step, for which we proposed an exact and a greedy strategy. We compared our framework with state-of-the-art local and global methodologies, showing its effectiveness and efficiency on well-known and synthetic datasets. Possible future research directions include: mechanisms in the greedy assignment step so as to maximize the chances of final \textit{feasible} assignments; \texttt{S-MDEClust} extensions for other clustering models; management of constraints different from the instance-level ones.

\section*{Competing Interest}
The authors have no competing interests to declare.

\section*{Data and Code Availability}
The \texttt{S-MDEClust} source code, the synthetic datasets, the constraint configurations for the $2^{\text{nd}}$ and $3^{\text{rd}}$ dataset sets (Table \ref{tab::datasets}) and all the computational experiments results can be found at \href{https://github.com/pierlumanzu/s_mdeclust}{github.com/pierlumanzu/s\_mdeclust.}




\bibliographystyle{elsarticle-num}

\end{document}